\documentclass[preprint]{elsarticle}
\usepackage{amssymb,amsfonts}
\usepackage{natbib}
\usepackage{amsthm}
\usepackage{amsmath}
\usepackage{graphicx}
\usepackage{caption}
\usepackage{subcaption}
\usepackage[active]{srcltx}
\usepackage{multirow}
\usepackage{hyperref}

 \newtheorem{theorem}{Theorem}[section]

 \theoremstyle{remark}
 \newtheorem{remark}[theorem]{\bf{Remark}}

\newcommand{\kbs}{k^*}
\newcommand{\Kbs}{K^*}
\newcommand{\Kb}{K}
\newcommand{\ks}{k^*}
\newcommand{\Phis}{\Phi^*}

\newcommand{\gs}{G^*}
\newcommand{\ps}{p^*}

\newcommand{\ubs}{\ub^*}

\newcommand{\us}{u^*}

\newcommand{\kbar}{\overline{K}}

\newcommand{\e}{\varepsilon}

\newcommand{\xib}{\mbox{$\xi$}}
\newcommand{\kb}{{k}}
\newcommand{\ub}{{u}}

\newcommand{\xb}{{x}}
\newcommand{\Xb}{{X}}
\newcommand{\av}[1]{\left<{#1}\right>}
\newcommand*{\doi}[1]{\href{http://dx.doi.org/\detokenize{#1}}{doi: \detokenize{#1}}}

\begin{document}

 \begin{frontmatter}
%
%
%
 \title{Upscaling of  Nonlinear Forchheimer Flows}

 \author{Eugenio Aulisa}
 \ead{eugenio.aulisa@ttu.edu}
 \address{Texas Tech University, Department of Mathematics and Statistics, Broadway and Boston, Lubbock, TX 79409-1042}

\author{Lidia Bloshanskaya\corref{cor1}}
\ead{lidia.bloshanskaya@ttu.edu}
\address{Texas Tech University, Department of Mathematics and Statistics, Broadway and Boston, Lubbock, TX 79409-1042}
\author{Yalchin Efendiev}
\ead{efendiev@math.tamu.edu}
\address{Texas A\&M University, Department of Mathematics,
College Station, TX 77843-3368}
\author{Akif Ibragimov}
\ead{akif.ibraguimov@ttu.edu}
\address{Texas Tech University, Department of Mathematics and Statistics, Broadway and Boston, Lubbock, TX 79409-1042}
\cortext[cor1]{Corresponding author}

%
%
%
%
\begin{abstract}
In this work we propose upscaling method for nonlinear 
Forchheimer flow in highly heterogeneous porous media.
The generalized Forchheimer law is considered for incompressible and 
slightly-compressible single-phase flows.
We use recently developed analytical results \cite{ABHI1} and write
the resulting system in terms of a degenerate nonlinear flow 
equation for the pressure with the nonlinearity that depends
on the pressure gradient. 
The coarse scale parameters for the steady state problem  are 
determined so that the volumetric average of velocity of the flow 
in the domain on fine scale and on coarse scale are close enough.
A flow-based coarsening approach is used, where the equivalent permeability 
tensor is first evaluated following the streamline of the existing 
linear cases, and successively modified in order to take 
into account the nonlinear effects.
Compared to previous works 
\cite{pes-trykozko, dur-fard}, our approach 
relies on recent analytical results of Aulisa et al. \cite{ABHI1}
 and combines it with rigorous mathematical
upscaling theory for monotone operators.
The developed upscaling algorithm for nonlinear steady state problems 
is effectively used for variety of heterogeneities in the
 domain of computation. 
Direct numerical computations for average velocity and productivity 
index justify the usage of the coarse scale parameters obtained for 
the special steady state case in the fully transient problem.
Analytical upscaling formulas in stratified domain are obtained for 
the nonlinear case. They correlate with high accuracy with numerical results.
\end{abstract}

 \end{frontmatter}

\section{Introduction}

In recent years, 
using 
near well data, e.g., core data, engineers create
increasingly complex and detailed
geocellular models which result
in highly heterogeneous geological parameters of reservoirs.
Such descriptions typically require a high number of computational cells
which are difficult to solve, e.g., in well optimization problems 
and history matching.
To reduce the computational complexity, some type of coarsening and
upscaling procedures are needed. 
The geological parameters, such 
as permeability or transmissibility and porosity, should be upscaled for each 
coarse-grid block. 

The variety of approaches for upscaling and multiscale methods of fine scaled geological parameters have been proposed for the linear Darcy case 
(e.g., \cite{hw97, eh09,  hughes98, jennylt03,
   aarnes04,     dur91,  weh02}). These approaches include
upscaling methods, e.g., \cite{dur91, weh02, eh09} and multiscale methods 
\cite{jennylt03, aarnes04, eh09}.
In both approaches, a goal is to represent the solution on a coarse grid
where each coarse-grid block consists of a union of connected fine-grid blocks.
In upscaling methods, the upscaled
permeability is calculated in each coarse-grid block
by solving local problems with specified boundary conditions
and calculating the average of the flow flux. 
Local problems can be solved
in extended domains for computing the effective properties.
In multiscale methods, the local multiscale basis functions are computed
instead of local effective properties and these basis functions
are coupled via a global formulation. 

The extensions of these methods to nonlinear flows, such as Forchheimer flow,
are carried out in several papers,  see \cite{dur-fard,pes-garib}
which are closely related to our work. In \cite{dur-fard}, the authors
consider the use of iterative upscaling techniques where at each
iteration, local-global upscaling technique is used. In \cite{pes-garib},
the authors use special nonlinear forms for upscaled Forchheimer flows
that simplify the upscaling calculations. In the current paper, our
goal is to carry out rigorous nonlinear upscaling using new formulations
of Forchheimer flows.

In current paper, we utilize recent finding \cite{ABHI1}, 
where Forchheimer
equation is written in an equivalent form using monotone nonlinear
permeability function depending on gradient of pressure. This equivalent formulation reduces the original system of equations for pressure and velocity to one nonlinear parabolic or elliptic equation for pressure only. The ellipticity constant of this equation degenerates as the pressure gradient converges to infinity. The rate of the degeneration is effectively controlled by the order of Forchheimer polynomial    and the structure of the coefficients has the important monotonicity properties \cite{ABHI1}, Proposition III.6 and Lemma III.10.
It allows to prove results on the well-posedness of the initial boundary value problem 
and apply numerical homogenization theory. 

In this paper we present the upscaling algorithm for fluid flow in incompressible media for two types of fluids, incompressible and slightly compressible.
Steady state problem for incompressible flow reduces to the degenerate elliptic equation, however the corresponding problem for compressible fluid reduces to time dependent degenerate parabolic equation.

In this report we first introduce and investigate the upscaling procedure for the time independent problem in case of incompressible fluid. 
In case of time dependent problem the question one should address is that while the solution is time dependent, the upscaled parameters are time independent for incompressible media.
We use the upscaled parameters obtained for steady state case  in  the
time dependent problem. 
This procedure is justified by the results obtained in our papers \cite{ABI11} and \cite{ABI12} and the numerical experiment presented in this article. 
Namely, we will relate the fine scale fully transient solution to the special pseudo steady state (PSS) solution. This solution has a form $At+W(x)$, where $A$ is a constant and $W(x)$ is a solution of auxiliary steady state boundary value problem for the equation with non zero RHS. 
According to our results in \cite{ABI11} and \cite{ABI12}
 under some assumptions 
the pseudo steady state pressure and velocity serve as pseudo attractors for fully transient pressure and velocity.
To upscale the steady state equation we determine the coarse scale porosity and nonlinear permeability, so that the average volumetric velocity of the flow is preserved.\\
To evaluate the described method for the time dependent case we compare the productivity index (PI) of the well on the fine  and  coarse grids. The PI is inversely proportional to the difference between the average of pressure in the reservoir and on the well. We select the PI as a criteria for the evaluation of the upscaling method  as it is widely used by the engineers \cite{Muskat, Raghavan}.
In the numerical examples we calculate the difference between the values of the PIs on fine and coarse grids. 
Our numerical results show that 
the proposed algorithm provides accurate results 
for different heterogeneities and nonlinearities in steady state case. 
Resulting transient velocity and PI on coarse scale also provide accurate approximation of corresponding transient parameters on fine scale
for heterogeneous fields considered in the paper. We expect the accuracy
of the proposed method depends on heterogeneities as in
a single-phase upscaling (\cite{cdgw03}), 
i.e., for highly heterogeneous
fields, the accuracy of the method will deteriorate. In this paper,
our main goal is to propose a method to handle the nonlinearities
and, thus, we do not consider very highly heterogeneous fields
(\cite{cdgw03}).\\
The paper is organized as follows.
In Sec.~\ref{sec:genforch} we introduce $g$-Forchheimer equations, review their properties and formulate the problem. 
In Sec.~\ref{sec:coarse-eq} we obtain the form of the coarse scale equation for generalized Forchheimer flow. Sec.~\ref{sec:coarse-compress} presents the discussion of convergence results for the transient velocity and PI in case of slightly compressible flow and usage of the upscaled parameters from the steady state equation in transient case.
Sec.~\ref{sec:algorithm} is devoted to description of upscaling algorithm. 
In Sec.~\ref{analytic} we obtain the explicit analytical upscaling formulas  in case of incompressible fluid for stratified region.
In Sec.~\ref{sec:num-incompres} and \ref{sec:num-compres} we present the numerical results for the incompressible and slightly compressible flows correspondingly.

\section{Problem statement and Preliminary results}\label{sec:problem-state}
\subsection{Generalized Forchheimer equation}\label{sec:genforch}
Darcy equation  describes the linear dependence of velocity $\ub$ on the pressure gradient $\nabla p$
\begin{equation}\label{darcy}
 \ub=-\tfrac1{\mu}\kb(x)\nabla p.
\end{equation}
Here  $\kb(x)$ is symmetric positive definite permeability tensor,
 $\mu$ is the viscosity of the fluid.

Forchheimer equation \cite{Forch} is known to generalize Darcy's equation to take into account inertial terms
and has been introduced in
the literature in several forms. E.g.,
\begin{equation}\label{Forch123}
\begin{aligned}
 &\text{Two term law:} && \ub+\beta(x)\|\ub\| \ub=-\tfrac1{\mu}\kb(x)\nabla p,\\
 &\text{Three term law:}&&\ub+a_1(x)\|\ub\| \ub +a_2(x) \|\ub\|^2 \ub=-\tfrac1{\mu}\kb(x)\nabla p,\\
 &\text{Power law:} &&\ub+b_1(x) \|\ub\|^{m-1} \ub=-\tfrac1{\mu}\kb(x)\nabla p, \qquad 1.6\leq m\leq2.
\end{aligned} 
\end{equation}
Coefficients $\beta(x)$, $a_1(x)$, $a_2(x)$ and $b_1(x)$ are empirical. 

All these relations can be written in a compact form as
\begin{equation}\label{g-forch}
 g(\|\ub\|,x)\ub=-\tfrac1\mu k(x)\nabla p,
\end{equation}
for some function $g(s,x)\geq0$ for $s\geq 0$. We will refer to \eqref{g-forch} as $g${\it-Forchheimer $($momentum$)$ equation}.
For simplicity from now on we assume the viscosity $\mu=1$, i.e. $\frac{\kb}{\mu}=\kb$.

To develop rigorous numerical homogenization concepts for Forchheimer flow,
we use the results in  \cite{ABHI1} which allows writing (\ref{g-forch}) as
a monotone relation for $\nabla p$. Moreover, this allows obtaining the well-posedness results of the corresponding initial boundary value problem 
 and allows estimating the residual error in numerical homogenization
because of monotonicity.
It was shown in \cite{ABHI1} that the monotone relation between velocity and
gradient of pressure exists for general functions  $g(s,x)$ in the form
\begin{equation}\label{gdef}
 g(s,x)=1+\sum_{j=1}^k a_j(x) s^{\alpha_j}=1 + a_1(x)s^{\alpha_1}+a_2(x)s^{\alpha_2}+\ldots +a_k(x)s^{\alpha_k},
\end{equation}
where $k\ge 0$,  the exponents satisfy $0<\alpha_j<\alpha_{j+1}$, and the coefficients $a_j(x)\ge0$, $j=1,\dots, k$. 
Thus defined function $g$ in \eqref{g-forch} includes all the known cases of Forchheimer flow \eqref{darcy} and \eqref{Forch123}.

We can define the inverse function
\begin{equation}\label{Kg}
 G(\xi;x)=\frac1{g(h^{-1}(\xi),x)},\quad \xi\geq0,\quad  h(s)=sg(s,x), \quad s\geq 0.
\end{equation}
We then can obtain the equivalent form of Eq.~\eqref{g-forch}
\begin{equation}\label{u-gdarcy-fine}
 \ub=-G(\|\kb(x)\nabla p\|;x)\,\kb(x)\,\nabla p,
\end{equation}
which we call {\it generalized $($nonlinear$)$ Darcy equation}.

%

\begin{remark}\label{remark-2forch}
 In the particular case of two-term Forchheimer law,
the nonlinear permeability tensor $G$  can explicitly be written
\begin{equation}\label{G-explicit}
 G(\xi; x)=\frac{2}{1+\sqrt{1+4 \beta\xi}}.
\end{equation}
\end{remark}

The $g$-Forchheimer equation written in the form \eqref{u-gdarcy-fine} allows reducing the dynamical system to single nonlinear equation of pressure. 
Namely, we  consider
the continuity equation
\begin{equation}\label{eq-cont}
 \phi(x)\frac{\partial\rho}{\partial t}=-\nabla\cdot(\rho \ub),
\end{equation}
where $\rho$ is the density of the fluid, and $\phi$ is the rock porosity.
For incompressible fluid ($\rho=const$), \eqref{eq-cont} reduces to $\nabla \cdot u=0$ and combined with the flow equation \eqref{u-gdarcy-fine} results in 
the degenerate elliptic equation of pressure only for steady-state flow
\begin{equation}\label{eq-incomp-fine}
 \nabla \cdot (G(\|\kb(x)\,\nabla p\|; x)\, \kb(x)\, \nabla p)=0.
\end{equation}

For slightly compressible fluid (such as the compressible liquid) the equation of state takes the form, see \cite{Muskat},
\begin{equation}\label{eq-state-compr}
\rho(p)=\rho_0 e^{\gamma p},
\end{equation}
where  $\gamma$ is the inverse of the compressibility constant.

Substituting \eqref{eq-state-compr} in \eqref{eq-cont} we get
\begin{equation}\label{eq-full}
 \phi \frac{\partial p}{\partial t}=-\tfrac1{\gamma}\nabla\cdot u +u \nabla p.
\end{equation}
 For slightly compressible fluids $\gamma$ is of order $10^{-8}$, thus 
we drop the second term in RHS of equation \eqref{eq-full}. Combining it with \eqref{u-gdarcy-fine} we obtain the degenerate parabolic equation for pressure 
\begin{equation}\label{degPar}
 \gamma \phi(x) \frac{\partial p}{\partial t}=\nabla \cdot  (G(\|\kb(x)\,\nabla p\|; x)\,\kb(x)\,\nabla p).
\end{equation}




Equations \eqref{eq-incomp-fine} and \eqref{degPar} 
describe the fluid flow on the fine grid.
 Our aim is to devise an
 upscaling algorithm for the parameters $k(x)$, $\phi(x)$ and $G(\|k\nabla p\|;x)$ in equations \eqref{eq-incomp-fine} and \eqref{degPar} and 
obtain the corresponding coarse scale equations.

%

\subsection{Coarse scale equation in case of incompressible fluid }\label{sec:coarse-eq}

To obtain the coarse scale equation for incompressible case,
 we first  rewrite Eq.~\eqref{eq-incomp-fine} in each coarse block $\Omega_c$ in a form:
\begin{equation}\label{eq-fin-gen}
 \nabla\cdot \Kb(\kb \nabla p;x)=0\ \ \text{in} \ \Omega_c,
\end{equation}
where
$\Kb(\eta;x)=G(\|\eta\|;x) \eta$.
We assume $p=\xi\cdot x$ on $\partial \Omega_c$, where vector $\xi=(\xi_1,\xi_2)$ (and $\xi=(\xi_1,\xi_2,\xi_3)$ in 3D).
We then solve Eq.~\eqref{eq-fin-gen} in each coarse block.
In each coarse block $\Omega_c$ we define
\begin{equation*}
 \Kbs( \xi)=\av{\Kb(\xi;x)}=\av{G(\|\xi\|; x)\, \xi}.
\end{equation*}
Here 
\begin{equation*}
\av{f}=\frac1{|\Omega_c|}\int_{\Omega_c}f\,d\Omega_c
\end{equation*}
is the volumetric average of the function over $\Omega_c$.

We would like to find the upscaled tensor $\kbs$ and scalar $\gs$, 
depending on $\kbs$, so that
\begin{equation*}
 \Kbs(\xi)=\gs(\xi)\kbs\xi.
\end{equation*}
Then, the upscaled equation takes the form
\begin{equation}\label{eq-incomp-coarse}
\nabla \cdot (\gs( \nabla \ps)\kbs \nabla \ps)=0.
\end{equation}

It follows that the coarse scale function $\gs$ depends on the vector $\kbs\nabla\ps$, while the fine scale function $G$ depends on the scalar $\|\kb(x)\,\nabla p\|$.

\subsection{Coarse scale equation in case of slightly compressible fluid}\label{sec:coarse-compress}

Unlike the steady state Eq.~\eqref{eq-incomp-fine} for incompressible fluid,  Eq.~\eqref{degPar} for slightly compressible flow is transient in time. 
Parameters $\kb$, $G$ and $\phi$ on the fine scale are, however, time independent. 
We want to find the upscaled parameters $\kbs$, $\gs$ and $\Phis$ on the coarse scale which are time independent as well. It is difficult to  use the original equation \eqref{degPar} for the upscaling  procedure  directly. Instead we will relate the fine scale transient pressure and velocity to the special pseudo steady state  solution of Eq.~\eqref{degPar} which will be defined below.

Let $U$ be the  domain with  the boundary $\Gamma$ consisting of two parts $\Gamma=\Gamma_e\cup\Gamma_i$.
The no-flux condition is imposed on $\Gamma_e$
\begin{equation}\label{bc-noflux}
\left.u\cdot\nu\right|_{\Gamma_e}=0; 
\end{equation}
and prescribed total flux condition is imposed on $\Gamma_i$ 
\begin{equation}\label{bc-totalflux}
\int_{\Gamma_i} \ub \cdot \nu\,ds = Q(t), 
\end{equation}
where $\ub$ is the velocity as in \eqref{u-gdarcy-fine} and $\nu$ is the outer normal to the boundary $\Gamma$.

In \cite{ABHI1} it was proved that there exists  a special solution $p_s(x,t)$ of equation \eqref{degPar} with  boundary condition
\eqref{bc-noflux} on $\Gamma_e$ such that 
\begin{equation}\label{pss-def}
 \frac{\partial p_s}{\partial t}=const.=-A \qquad\text{for all}\quad t.
\end{equation}
Such solution is  called \textit{Pseudo Steady State $($PSS$)$}.
From definition \eqref{pss-def}  of PSS solution it follows that the corresponding production rate is constant $Q(t)=Q=A|U|=const.$ and
 $p_s(x,t)$ can be written as
\begin{equation*}
 p_s(x,t)=-\gamma\,\frac{Q}{|U|} t + W(x) + C,
\end{equation*}
where $W(x)$ is called a \textit{basic profile} and is a solution of the steady state BVP
\begin{align}
& \nabla \cdot (G(\|\kb(x)\nabla W\|; x)\kb(x)\nabla W) =-\gamma\,\frac{Q}{|U|}\phi(x) ,\label{eq-compr-fine}\\
& \left. u_s(x)\cdot\nu\right|_{\Gamma_e}=0, \label{pss-bc-neu}\\
& \left. W\right|_{\Gamma_i}=\varphi_0(x), \label{pss-bc-dir}
\end{align}
with given function $\varphi_0(x)$ and constant $C$. 
Notice that $\nabla p_s=\nabla W$ and the corresponding PSS velocity
\begin{equation}\label{pss-vel}
 \ub_s(x)=-G(\|\kb\nabla p_s\|;x)\,\kb\,\nabla p_s=-G(\|\kb\nabla W\|;x)\,\kb\,\nabla W
\end{equation}
is time independent.

The steady state BVP \eqref{eq-compr-fine}-\eqref{pss-bc-dir} will be used  to find the upscaling parameters for fully transient equation \eqref{degPar}.
On coarse scale the steady state Eq.~\eqref{eq-compr-fine} will take the form
\begin{equation}\label{eq-compr-coarse}
 \nabla \cdot (\gs\kbs\nabla \ps) =-\frac{Q}{|U|}\,\Phis. 
\end{equation}

The upscaling algorithm for  $\kbs$, $\gs$ and $\Phis$ 
follows, as previously, a procedure
of equating average velocities (cf. \cite{dur91}). 
Consequently, 
\begin{equation}\label{uerr-abs}
 \|\av{\ub_s}_U-\av{\ubs_s}_U\|\quad\text{is sufficiently small.}
\end{equation}
Here 
\begin{equation}\label{av-vel}
\av{\ub_s}_U=\frac1{|U|}\int_{U}\ub_s \,dU; \qquad \av{\ubs_s}_U=\frac1{|U|}\int_{U} \ubs_s \,dU,
\end{equation}
where $\ub_s$ and $\ubs_s$ are the steady state velocities on fine and coarse scales, correspondingly.

If initial data is not of basic profile then solution $p(x,t)$ of the original equation \eqref{degPar} and the corresponding velocity $\ub(x,t)$ are time dependent. Thus in order to justify the upscaling criteria \eqref{uerr-abs} for general case 
 one should prove convergence of the corresponding time dependent quantity to the time independent one. This property was obtained in \cite{ABI11} and \cite{ABI12} under certain conditions on the boundary data.
Namely, let 
\begin{equation*}
 \psi(x,t)=p(x,t)|_{\Gamma_i}-\frac1{|\Gamma_i|}\int_{\Gamma_i}p(x,t)\,ds
\end{equation*}
and
\begin{equation*}
 \varphi(x)=\varphi_0(x)-\frac1{|\Gamma_i|}\int_{\Gamma_i}\varphi_0(x)\,ds 
\end{equation*}
be the deviations from the average on the boundary $\Gamma_i$ of the trace of transient solution $p(x,t)$ and basic profile $W(x)$ correspondingly.

We proved that if the differences $Q(t)-Q$ and $\psi(x,t)-\varphi(x)$ converge in certain sense to zero at time infinity  (see \cite{ABI12}, \S 3.2), then the PSS  velocity $\ub_s(x)$ serves 
as the steady-state attractor for the fully transient velocity  $\ub(x,t)$ with any initial data: 
\begin{equation*}
 \int_U\|\ub(x,t)- \ub_s(x)\|^2\,dx\to 0\qquad \text{as}\quad t\to\infty. 
\end{equation*}
This justifies the usage of criteria \eqref{uerr-abs} for the upscaling of coefficients $\kb$, $G$, and $\phi$  in fully transient problem.

To evaluate this method we made a comparison between  the {\it productivity index}  on coarse and fine scale using the  coefficients $\kbs$, $\gs$, and $\Phis$ on coarse scale.
Productivity index is routinely used by engineers in estimation of available reserves and optimizing well recovery efficiency (see \cite{Muskat, Raghavan, Slider}). It is defined as follows.
Let $p(x,t)$ be the solution of BVP in region $U$ for equation \eqref{degPar} with boundary conditions \eqref{bc-noflux} and \eqref{bc-totalflux}.
The \textit{Productivity Index/Diffusive Capacity $($PI$)$} is defined  as the ratio
\begin{equation}\label{PI-def}
 J(t)=\frac{Q(t)}{\overline{p}_U(t)-\overline{p}_{\Gamma_i}(t)},
\end{equation}
where $\overline{p}_U(t)-\overline{p}_{\Gamma_i}(t)$ is a {\it pressure drawdown}; and
\begin{equation*}
\overline{p}_U(t)=\frac1{|U|}\int_Up(x,t)\,dx,\qquad \overline{p}_{\Gamma_i}(t)=\frac1{|\Gamma_i|}\int_{\Gamma_i}p(x,t)\,ds.
\end{equation*}
%
%
%
%
%

In case of PSS flow the productivity index $J(t)$ is time independent and 
\begin{equation}\label{pi-pss}
 J(t)=J_{PSS}=\frac{Q}{\frac{1}{|U|}\int_U W(x)\,dx}.
\end{equation}

We compare the $J_{PSS}$ corresponding to the solution of the fine scale equation \eqref{eq-compr-fine}  with the $J^*_{PSS}$ corresponding to the solution of the equation \eqref{eq-compr-coarse} on coarse scale . 
It is numerically proved that the difference $|J_{PSS}-J^*_{PSS}|$  is small (see Sec.~\ref{sec:num-compres}).   

As it has been already mentioned, in general, the productivity index is time dependent. As in the case for velocity, it was proved  that  if the differences $Q(t)-Q$ and $\psi(x,t)-\varphi(x)$ converge in certain sense to zero at time infinity  (see \cite{ABI12}, \S 3.2), then
\begin{equation*}
 |J(t)- J_{PSS}|\to0\qquad \text{as}\quad t\to\infty. 
\end{equation*}
Thus  the coarse coefficients $\kbs, \gs$ and $\Phis$, obtained for the steady state equation \eqref{eq-compr-fine}, can be used to calculate fully transient productivity index on coarse scale. 

\begin{figure}[!h]
\begin{center}
 \includegraphics[height=65mm]{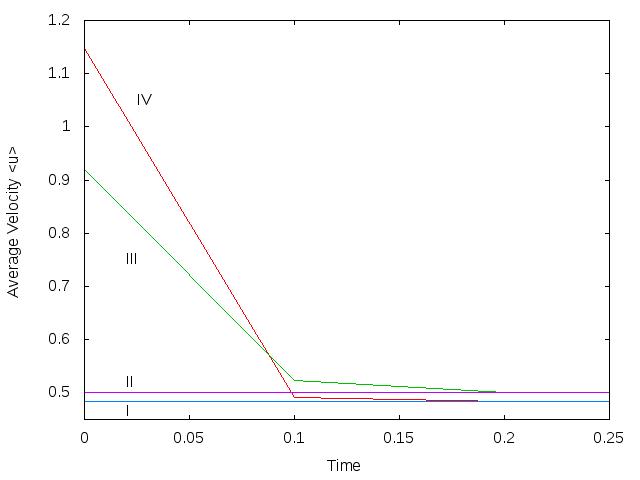}
\caption{Time dependence of the average velocity on fine and coarse scales: I - PSS $\left<u_s\right>$ on the fine scale; II - PSS $\left<u_s\right>$ on the coarse scale; III - $\left<u(t)\right>$ on the coarse scale; IV - $\left<u(t)\right>$ on the fine scale}
\label{fig:vel-timedep}
 \includegraphics[height=65mm]{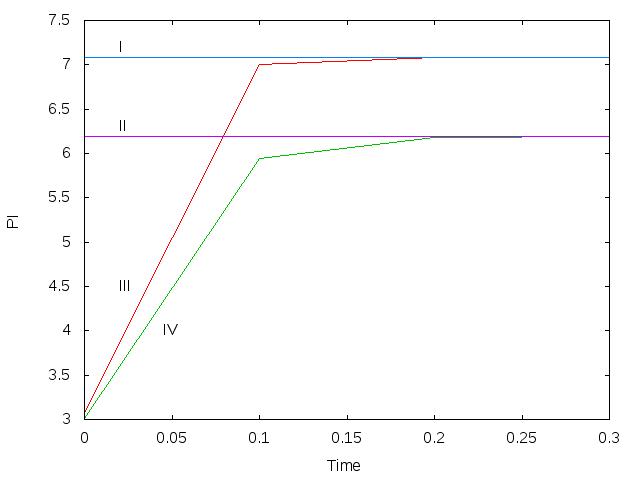}
\caption{Time dependence of the Productivity Index on the fine and coarse scales: I - PSS PI on the fine scale; II - PSS PI on the coarse scale; III - PI($t$) on the fine scale; IV - PI($t$) on the coarse scale}
\label{fig:pi-timedep}
\end{center}
\end{figure}

Numerical experiment confirms the theoretical findings above. Figures \ref{fig:vel-timedep} and \ref{fig:pi-timedep} present the time dependence of velocity and the PI of coarse and fine scales. The time dependent values are also compared to the PSS values, which are constant in time. As it can be seen from the graphs, in the long term the coarse scale time dependent velocity and PI calculated using the upscaled parameters from the steady state problem provide good approximation of the corresponding fine scale values.

\section{Numerical upscaling algorithm}\label{sec:algorithm}

\begin{figure}[ht]
\begin{center}
\includegraphics[height=50mm]{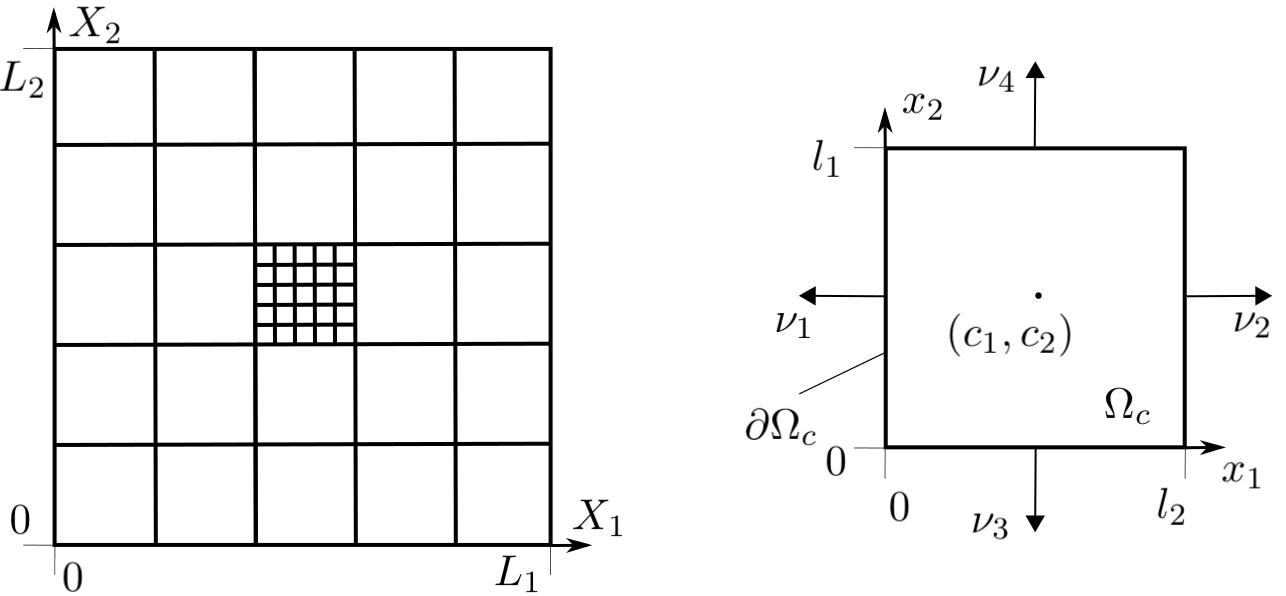}
\caption{Fine and Coarse Scale}
\label{fig:fine_and_coarse}
\end{center}
\end{figure}

In this section we present the numerical upscaling algorithm for the steady state equations \eqref{eq-compr-fine} and  \eqref{eq-incomp-fine}. We consider 2D rectangular region $\Omega$, with horizontal size $L_1$ and vertical size $L_2$ and two orthogonal grids: fine $\xb=(x_1,x_2)$-scale and coarse 
$\Xb=(X_1,X_2)$-scale (see Fig.~\ref{fig:fine_and_coarse}).

%

The porous media on the fine scale is considered  to be isotropic, and permeability tensor $\kb(\xb)$ is a  scalar function $k(x_1,x_2)$. 
The fine scaled equation \eqref{eq-compr-fine} or   \eqref{eq-incomp-fine} with parameters $k(x_1,x_2)$, $\phi(x_1,x_2)$ and
  $G(\|k\nabla p\|;x_1,x_2)$ is upscaled to the coarse scale equation \eqref{eq-compr-coarse} or  \eqref{eq-incomp-coarse} with parameters $\kbs(X_1,X_2)$, $\Phis(X_1,X_2)$ and
$\gs(\nabla p^*;X_1,X_2)$ so that condition \eqref{uerr-abs} is satisfied. The $u_s$ and $u_s^*$ in \eqref{uerr-abs} are the velocity on the fine and coarse scale correspondingly.


Our approach is purely local, so the algorithm is described for single coarse block $\Omega_c$ with boundary $\partial \Omega_c$.
For simplicity we take $\Omega_c$ to be the  rectangular
$[0,l_1]\times[0,l_2]$ with the area $|\Omega_c|=l_1\cdot l_2$, see Fig.~\ref{fig:fine_and_coarse}.

For each coarse block the two-step procedure is performed:
\vspace{-0.1cm}
\begin{itemize}
 \item[Step 1:]the equivalent permeability tensor $\kbs$ is obtained using linear upscaling methods;
 \item[Step 2:]the equivalent nonlinear coefficient $\gs$ is obtained using $\kbs$.
\end{itemize}

\textit{Step $1.$ Procedure to obtain  $\kbs$. } 
In order to obtain full permeability tensor $\kbs$ we use the standard local procedure via volume averages of velocity and pressure gradients, see for example \cite{Durlofsky}.
 We solve two flow problems in each block with periodic boundary conditions.
Namely, let $p_1$ and $p_2$ be the solutions of the fine scale equation in coarse block $\Omega_c$
\begin{equation}\label{eq_p_k}
 \nabla\cdot(k(x_1,x_2)\nabla p)=0,
\end{equation}
with boundary conditions:
\begin{equation}\label{k-bc}
  \left\{
\begin{aligned}
  & p_1(x_1,0)=p_1(x_1,l_2) \ \ \text{for}\ \ x_1\in[0,l_1];\\
 & \ub_1(x_1,0)\cdot \nu_3=-\ub_1(x_1,l_2)\cdot \nu_4;\\
 & p_1(0,x_2)=0;\quad p_1(l_1,x_2)=1;
\end{aligned}
\right. \ 
  \left\{
\begin{aligned}
 & p_2(0,x_2)=p_2(l_1,x_2) \ \ \text{for}\ \ x_2\in[0,l_2];\\
& \ub_2(0,x_2)\cdot \nu_1=-\ub_2(l_1,x_2)\cdot \nu_2;\\
 & p_2(x_1,0)=0;\quad p_2(x_1,l_2)=1.
\end{aligned}
\right.
\end{equation}
Here $\ub_i$ is the velocity vector corresponding to the pressure distribution $p_i$, $i=1,2$.

The four elements of the upscaled permeability $\kbs$ are then calculated from two vector equations: 
\begin{equation*}
 \av{\ub_i}=-\kbs\av{\nabla p_i}, \quad i=1,2.
\end{equation*}

The upscaled porosity $\Phis$ is computed via integral averaging on the coarse block following classical approach, e.g. \cite{Durlofsky}:
\begin{equation}\label{phia}
\Phis=\av{\phi(x_1,x_2)}=\frac1{|\Omega_c|}\int_{\Omega_c}\phi(x_1,x_2)\,d\Omega_c.
\end{equation}

\vspace{0.3cm}

\textit{Step $2:$ Procedure to obtain $\gs$. }
We use the upscaled permeability $\kbs$ to determine the nonlinear coefficient $\gs$
via pure local averaging.
As it was mentioned in Sec.~\ref{sec:coarse-eq}, unlike the fine-scale function $G$ depending on $\|k\nabla p\|$, the upscaled $\gs$ depends on the vector $\nabla \ps$ itself.
Let $\xib=(\xi_1,\xi_2)$ be the gradient of pressure in coarse block  $\Omega_c$.
For fixed $\xib$,  the $\gs$ is a constant.
If $\ps$ is the solution of coarse scale Eq.~\eqref{eq-incomp-coarse} with boundary condition
 \begin{equation}\label{xi-bc-lin}
\ps|_{\partial \Omega_c}=\xi_1 x_1+\xi_2 x_2,
\end{equation}
then it is also the solution of equation 
$  \nabla\cdot(\kbs\nabla \ps)=0$
with the same boundary condition. 

We determine $\gs$ so that 
\begin{equation}\label{gs-cond}
\|\av{\ubs}\|=\|\av{\ub}\|,
\end{equation}
where $\ubs=-\gs\kbs\nabla\ps$ is the velocity on coarse scale and $\ub=-G\kb\nabla p$ is the  velocity on fine scale corresponding to the solution $p(x,t)$ of   \eqref{eq-incomp-fine} with boundary condition
$p|_{\partial \Omega_c}=\xi_1 x_1+\xi_2 x_2$.
Then  for fixed $\xi_1, \xi_2$ we have:
\begin{equation}\label{gs}
 \gs(\xi_1,\xi_2)=\frac{\|\av{\ub}\|}{\|\kbs \av{\nabla \ps}\|}.
\end{equation}
 
Using formula \eqref{gs} we numerically construct the table of values of $\gs$ for $\xi_1,\xi_2\in(-\infty,\infty)$.
It follows that $\gs(0,0)=1$,
$\gs(\xi_1, \xi_2)\to 1$ if $\|\xib\|\to 0$, $\gs(\xi_1,\xi_2)\to 0$ if $\|\xib\|\to \infty$
and  $\gs$ possesses certain symmetry:
$ \gs(\xi_1,\xi_2)=\gs(-\xi_1,-\xi_2)$.
It is thus sufficient to consider $\xi_1\in(-\infty,\infty)$ and $\xi_2\geq 0$ only.
It is worth mentioning that the special attention should be paid to the way the domain for the $\xi$ is discretized.
Taking the grid to be too fine makes the calculations overly expensive, however the sparse grid does not allow to capture the features of nonlinearity of the process.
We will use the non uniform grid, where the subsequent point is calculated on the basis of the deviation between the preceding values of the function. Namely,  $\xi_1$ and $\xi_2$ are taken from the set $\eta_n$, $n=0,1,2,\dots$,

where the first three values are taken \textit{a priori}:
 $\eta_0=0$, and $\eta_1$, $\eta_2$ to be small enough. Next value $\eta_{n+1}$ is chosen so that 
\begin{equation*}
 \frac{\gs_n-\gs_{n+1}}{\gs_n}\leq\e_n,\quad\text{where}\quad \gs_m=\gs(\eta_m,\eta_m),\quad m=n,n+1
\end{equation*}
for some set value $\e_n$.
The stopping criteria for the computation is $\gs_n\le\e$ and $|(\gs_n)'|\le\e_d$ so that the value of $\gs$ as well as its variation are sufficiently small.

\begin{figure}[!h]
\begin{center}
\includegraphics[height=50mm]{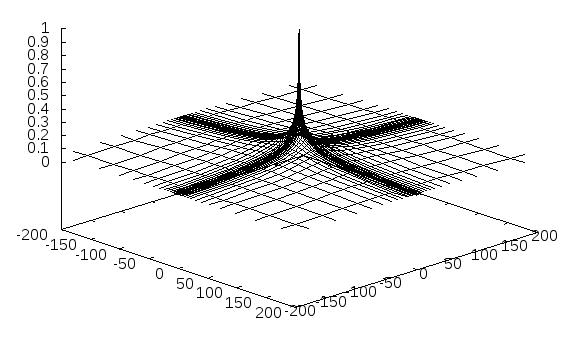} 
\caption{The shape of function $\gs(\xi_1,\xi_2)$}
\label{gs-func}
\end{center}
\end{figure}

The shape of function $\gs$ is presented on Fig.\ref{gs-func}.

\section{Analytical Upscaling for the Layered Porous Media}\label{analytic}

\begin{figure}[ht]
\begin{center}
\includegraphics[height=40mm]{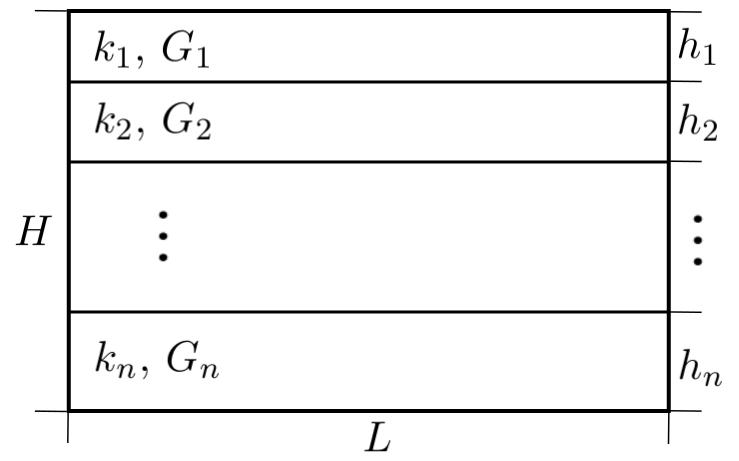}
\caption{Layered porous media, region $R$}
\label{fig_cell}
\end{center}
\end{figure}

Here we will present the analytical upscaling formula for the nonlinear Forchheimer flow of incompressible fluid in layered porous media.
Consider a rectangular region $R$ of horizontal size $L$ and vertical size $H$. 
The region has a horizontal multilayer structure, and is composed by $n$ layers, see Fig.~\ref{fig_cell}.
Each $i$th layer, $i=1,\dots, n$, has vertical size $h_i$ and is characterized 
by constant isotopic permeability $k_i$ and $g$-Forchheimer polynomial $g(s,x)=g_i(s)$ with constant coefficients or, equivalently, by the nonlinear function $G_i=G_i(k_i\|\nabla p_i\|)$. 
We assume that the type of nonlinearity is the same for each layer, while the coefficients of $g$-polynomials can be different.

Under these assumptions, in each $i$th layer, $i=1,\dots, n$, equations \eqref{g-forch}, \eqref{u-gdarcy-fine} and \eqref{eq-incomp-fine}
yield
\begin{align}
&g_i(\|u_i\|)\,u_i=- k_i\nabla p_i,\label{g-forch-i}   \\
&u_i=-G_i k_i \nabla p_i, \label{GDarcyi}\\
&\nabla \cdot (G_i k_i\nabla p_i)=0.\label{degElli}
\end{align}
Here $\ub_i$, $p_i$ and $Q_i$ are, correspondingly, velocity, pressure and the total boundary flux in $i$th layer $R_i$ and $\|\cdot\|$ is the $l_2$ vector norm.

We assume that flow within the whole block $R$ is subject to the equation with the same type of nonlinearity as in each layer.
We aim to find the equivalent homogeneous block permeability $\ks$ and nonlinear coefficient $\gs=\gs(\ks\nabla \ps)$ for two types of flow: flow parallel to the layers 
(Sec.~\ref{sec:incompress-horizontal}) and flows perpendicular to the layers (Sec.~\ref{sec:incompress-vertical}).
The upscaled parameters are determined so that the total flux of the system stays the same as with nonhomogeneous parameters.
The comparison between the obtained analytical results and numerical computations using the method in Sec.~\ref{sec:algorithm} are presented in Sec.~\ref{sec:num-incompres}.
Note, that the case when $g_i=const.$ is the same as Darcy case and the upscaling formulas for $\ks$ are the same as in \cite{Bear}.

\subsection{Flow Parallel to the Layers}\label{sec:incompress-horizontal}
We impose the following boundary conditions on boundaries of $R$
\begin{itemize}
\item $p_i=\ps=p_0$,  on the left boundary, $i=1,\dots,n$,
\item $p_i=\ps=p_1$ on the right boundary, $i=1,\dots,n$,
\item $\ub \cdot \nu =0 $ on the bottom and top boundaries,
\end{itemize}
where $\nu$ is the outward normal, and $p_1>p_0$.
Under these conditions the flow is parallel to the layers and the solution of Eq.~\eqref{degElli} is  linear in $x$.
The pressure gradient is constant and is equal to $\nabla p=(\xi,0)$, where $\xi=(p_1-p_0)/L$.
In each layer the vertical velocity component $u_{i_2}$ is identically zero,
while the horizontal component $u_{i_1}$ is constant in each layer and, according to \eqref{GDarcyi}, is equal
\begin{equation}\label{vi}
u_{i_1}=u_i= -G_i(k_i\xi)\, k_i\, \xi.
\end{equation}
On the other hand the outgoing flux is equal to incoming flux and is equal  to the sum of fluxes in the $i$th layer:
\begin{equation}\label{hor-flux}
 Q=\sum_{i=1}^n Q_i=- \sum_{i=1}^n  u_{i} h_i = \xi \sum_{i=1}^n G_i k_i h_i.
\end{equation}
The flux is zero on the top and bottom boundaries.

We now consider the analogous block with the same boundary conditions and  permeability $\ks$ and nonlinear function $\gs$ 
resulting in the same flux $Q$.
In this case the flux is 
\begin{equation*}
Q=-u^*H=\gs k^* H\, \xi,
\end{equation*}
where $u^*$ is constant horizontal component of upscaled velocity.
Expression for $Q$ above and \eqref{hor-flux} yield
\begin{equation}\label{EqtFlx}
\xi \sum_{i=1}^n G_i k_i h_i=\gs k^* H\, \xi\qquad\text{or}\qquad \us=\frac1H\sum_{i=1}^n u_ih_i.
\end{equation}
First we consider the limiting linear Darcy case $G_i=1$. In this case  $\gs=1$.
We then find an expression for $k^*$
\begin{equation}\label{||ks}
 k^*=\frac{1}{H}\sum_{i=1}^n {k_i h_i}.
\end{equation}
In view of \eqref{||ks} the general expression for $\gs$ follows from \eqref{EqtFlx}
\begin{equation}\label{||Gs}
 \gs=\frac1H\frac1{\ks}\sum_{i=1}^n G_i k_i h_i=\frac{\sum_{i=1}^n G_i k_i h_i}{\sum_{i=1}^n k_i h_i}.
\end{equation}
Formulas \eqref{||ks} and \eqref{||Gs} can be generalized in case when the parameters $k=k(x_2)$, $G=G(k\xi,x_2)$ are continuous functions:
\begin{equation*}
  \ks=\frac{1}{H}\int_{0}^H k(x_2)\,dx_2;\qquad \gs=\frac1H\frac1{\ks}\int_0^H G(k\xi,x_2)k(x_2)\,dx_2.
\end{equation*}

Alternatively, using $g$-Forchheimer equation \eqref{g-forch-i} with $\|u_i\|=u_i$, the upscaling formula for the $g$-polynomial can be obtained:
\begin{equation*}
\frac1{g^*(\us)}=\frac{\sum_{i=1}^n \frac{1}{g_i(u_i)} k_i h_i}{\sum_{i=1}^n k_i h_i}.
\end{equation*}

From here we can obtain the upscaled coefficients $a^*_j$, $j=1,\dots,m$ corresponding to the power $s^{\alpha_j}$, for the $g$-polynomial in domain $R$.
In particular in case of two-terms law as in Remark~\ref{remark-2forch} the upscaled Forchheimer coefficient can be obtained explicitly in the form
\begin{equation}\label{||betas-expl}
 \beta^*=
\frac{\sum_{i=1}^n  \beta_i u_{i}^2 \frac{h_i}{H}}{ \left(\sum_{i=1}^n u_{i} \frac{h_i}{H}\right)^2}
=\dfrac{\sum_{i=1}^n  \beta_i \left( \frac{2}{1+\sqrt{1+4 \beta_i k_i \xi } }k_i\right)^2 \frac{h_i}{H}}
{ \left(\sum_{i=1}^n \frac{2}{1+\sqrt{1+4 \beta_i k_i \xi } }k_i \; \frac{h_i}{H}\right)^2},
\end{equation}
where $\beta_i$ is coefficient corresponding to $i$th layer.

The coefficient $\beta^*$ depends explicitly on $\xi$.
The two limiting cases are
\begin{equation*}
\lim_{\xi \rightarrow 0}\beta^* =
\frac{\sum_{i=1}^n \beta_i k_i^2 \frac{h_i}{H} } { \left(\sum_{i=1}^n k_i h_i\right)^2 }=
\frac{\sum_{i=1}^n \beta_i k_i^2 \frac{h_i}{H} } { k^{*2} }
\end{equation*}
and
\begin{equation*}
\lim_{\xi \rightarrow \infty}\beta^* =
\frac{\sum_{i=1}^n k_i \frac{h_i}{H}}{\left( \sum_{i=1}^n \sqrt{\frac{k_i}{\beta_i}} \frac{h_i}{H} \right)^2}=
\frac{k^*}{\left( \sum_{i=1}^n \sqrt{\frac{k_i}{\beta_i}} \frac{h_i}{H} \right)^2}.
\end{equation*}

In case when the parameters  $k=k(x_2)$ and $\beta=\beta(x_2)$ are continuous functions, the expression
\eqref{||betas-expl} for $\beta^*$ yields
\begin{equation*}
\beta^*= 
\dfrac{\frac{1}{H}\int_{0}^H  \beta(x_2) \left( \frac{2 k(x_2) }{1+\sqrt{1+4 \beta(x_2) k(x_2) \xi } }    \right)^2\,dx_2}
{ \left(\frac{1}{H}\int_{0}^H \frac{2 k(x_2) }{1+\sqrt{1+4 \beta(x_2) k(x_2) \xi } }\,dx_2\right)^2}.
\end{equation*}

\subsection{Flow Perpendicular to the Layers}\label{sec:incompress-vertical}
Let consider the same geometry and let impose the following boundary conditions
\begin{itemize}
\item $p_i|_{x_2=0}=\ps|_{x_2=0}=p_0$, on the bottom boundary, 
\item $p_i|_{x_2=H}=\ps|_{x_2=H}=p_n$, on the top boundary,
\item $\ub \cdot \nu =0 $ on the left and right boundaries.
\end{itemize}

In this case the flow is perpendicular to the layers and the horizontal velocity component $u_{i_1}$ is identically zero, while
the  vertical component of velocity $u_{i_2}=u$ is constant in each layer. 
The pressure gradient in $i$th layer is equal to 
\begin{equation*}
 \nabla p_i=(0,\xi_i),\quad\text{where}\quad  \xi_i=\frac{p_i-p_{i-1}}{h_i},
\end{equation*}
where $p_i$ is the pressure measured at the top of $i$th layer for $i=1,\dots,n$.
Then, according to \eqref{GDarcyi}, vertical component of velocity is equal to
$ u=-G_i k_i \xi_i$.

It thus follows that the flux is constant and in each layer is equal to  $Q=-uL=G_i k_i \xi_i L$.
Then the pressure gradient in $i$th layer is
\begin{equation}\label{xi-i}
 \xi_i=\frac{Q}{G_ik_i L}, \quad i=1,\dots,n.
\end{equation}

We again want to identify the equivalent homogeneous parameters $\ks$ and $\gs$ in the region $R$ resulting in the same flux $Q^*=Q$.
The pressure gradient in the domain $R$ is 
\begin{equation}\label{sxi}
\xi= \frac{p_n-p_0}H= \frac{1}{H}\sum_{i=1}^n\xi_i h_i. 
\end{equation}
We get the expression for the flux
\begin{equation}\label{qs}
 Q=-u^*L=\gs \ks \xi L.
\end{equation}
Plugging  \eqref{sxi} in \eqref{qs} and using \eqref{xi-i}  we get
\begin{equation*}
 \frac1{\gs\ks}=\frac1H\sum_{i=1}^n\frac1{G_i}\frac{h_i}{k_i}
\end{equation*}
and it follows:
\begin{equation}\label{|ks-gs}
  \frac{1}{\ks}=  \frac{1}{H}\sum_{i=1}^n\frac{h_i}{k_i};\qquad\qquad  \frac1{\gs}=\frac{\ks}{H}\sum_{i=1}^n \frac1{G_i}\frac{h_i}{k_i}.
\end{equation}

Alternatively, using $g$-Forchheimer equation \eqref{g-forch-i} $g_i(u)u=-k_i\xi_i$, the upscaling formula for the $g$-polynomial can be obtained:
\begin{equation*}
 g^*(u)= \frac{\ks}{H}\sum_{i=1}^ng_i(u)\frac{h_i}{k_i}=\frac{\sum_{i=1}^ng_i(u)\frac{h_i}{k_i}}{\sum_{i=1}^n\frac{h_i}{k_i}},
\end{equation*}
and thus 
\begin{equation}\label{|ajs}
 a^*_j= \frac{\ks}{H}\sum_{i=1}^na_{j,i}\frac{h_i}{k_i}=\frac{\sum_{i=1}^na_{j,i}\frac{h_i}{k_i}}{\sum_{i=1}^n\frac{h_i}{k_i}}, \quad j=1,\dots, m.
\end{equation}
where $a_{j,i}$ is the coefficient of $g$-polynomial corresponding to power $s^{\alpha_j}$, $j=1,\dots, m$ (see Eq.~\eqref{gdef}) in $i$th layer.

In  case when the parameters  $k=k(x_2)$, $G=G(k \xi(x_2),x_2)$ and $a_{j}=a_j(x_2)$, $j=1,\dots, m$, are continuous functions in $x_2$, Eqs.~\eqref{|ks-gs} and \eqref{|ajs}
yield
\begin{equation*}
\frac{1}{k^*}=\frac{1}{H} \int_{0}^H \frac{dx_2}{k(x_2)};\ 
\frac1{\gs}=\frac{\ks}{H}\int_0^H\frac{dx_2}{G(k\xi(x_2),x_2)k(x_2)};\ 
 a^*_j= \frac{\ks}{H}\int_0^H \frac{a_{j}(x_2)}{k(x_2)}\,dx_2.
\end{equation*}
%

\section{Numerical Results}

\begin{figure}
 \centering
\begin{subfigure}[b]{0.49\textwidth}
                \centering
                \includegraphics[height=40mm]{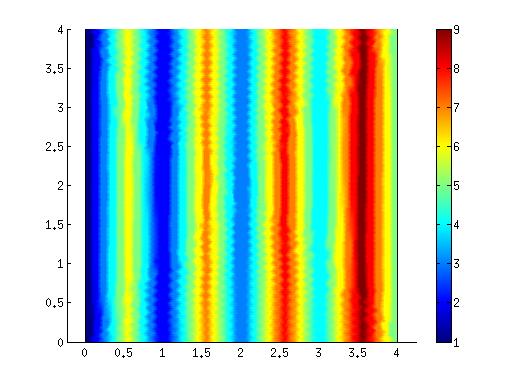}
                \caption{Vertically stratified $k(x_1,x_2)$}
                \label{fig:kver}
\end{subfigure}
\begin{subfigure}[b]{0.49\textwidth}
                \centering
                \includegraphics[height=40mm]{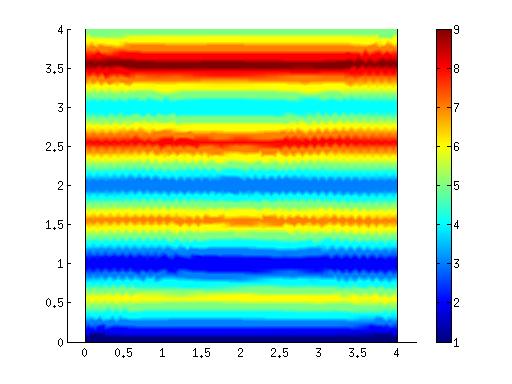}
                \caption{Horizontally stratified $k(x_1,x_2)$}
                \label{fig:khor}
\end{subfigure}

\begin{subfigure}[b]{0.49\textwidth}
                \centering
                \includegraphics[height=40mm]{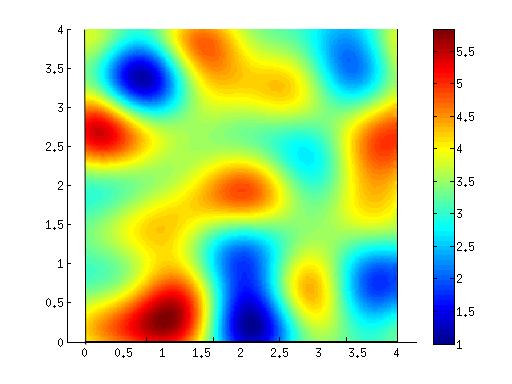}
                \caption{Randomly distributed $k(x_1,x_2)$}
                \label{fig:krand}
\end{subfigure}
\begin{subfigure}[b]{0.49\textwidth}
                \centering
                \includegraphics[height=40mm]{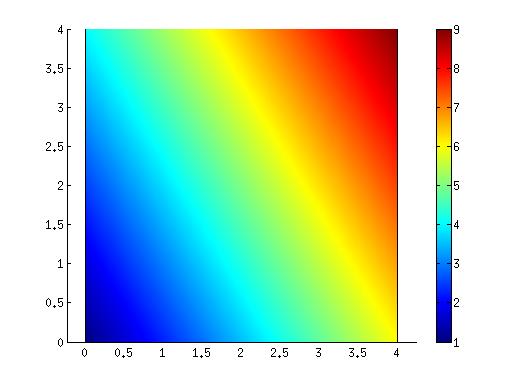}
                \caption{Linearly distributed $k(x_1,x_2)$}
                \label{fig:kplane}
\end{subfigure}
\caption{Permeability $k$ on fine scale}
\label{fig:k-distr}
\end{figure}

In this section we numerically illustrate the upscaling algorithm described in Sec.~\ref{sec:algorithm} for the incompressible and slightly compressible fluids. 
The considered cases of the permeability distribution on fine scale are presented in Fig.~\ref{fig:k-distr}.
The obtained upscaling errors are relatively small, since we did not consider large heterogeneities, but instead focused on the upscaling method for the nonlinear flow.
\subsection{Numerical Results for Incompressible Fluid}\label{sec:num-incompres}
In this section we present the numerical results for  incompressible flow.
Several approaches are compared: the upscaling algorithm Sec.~\ref{sec:algorithm} and the analytical formulas obtained in Sec.~\ref{analytic}.

On the fine scale the pressure is subject to equation \eqref{eq-incomp-fine} in the region $\Omega=[0,L_1]\times[0,L_2]$ with the boundary 
conditions
\begin{equation*}
 p(0,x_2)=0; \quad p(L_1,x_2)=1;\quad \tfrac{\partial p}{\partial x_2}(x_1,0)=\tfrac{\partial p}{\partial x_2}(x_1,L_2)=0.
\end{equation*}

We report the relative error in the averaged velocities originated from the upscaling of equation \eqref{eq-incomp-fine} to \eqref{eq-incomp-coarse}:
\begin{equation}\label{vel-err}
\frac{\|\left<\ub\right>-\left<\ubs\right>\|}{\|\left<\ub\right>\|}.
\end{equation}
The errors for the layered system are reported in Table~\ref{tab:incomp-||flow} (flow parallel to the layers, permeability as in Fig.~\ref{fig:khor}) and Table~\ref{tab:incomp-|flow} (flow perpendicular to the layers, permeability as in Fig.~\ref{fig:kver}).
The results for the case of random system are reported in Table~\ref{tab:incomp-rand}.

For each case we compare results obtained in three different ways: 
1) analytical formulas Eqs.~\eqref{||ks}, \eqref{||Gs} for  flow parallel 
to the system (denoted by ``Av $\Rightarrow$'' in the tables);
2) analytical formulas Eqs.~\eqref{|ks-gs} 
for  flow perpendicular to the system (denoted by ``Av $\Uparrow$'' in the tables);
3) numerical approach described in Sec.~\ref{sec:algorithm} (denoted by ``Num'' in the tables). 
For the layered system,  
the corresponding analytical formula gives the exact result.

The calculations are performed for different orders of nonlinearity:\\
\indent 1. linear Darcy case, $\beta =0$;\\
\indent 2. two-term Forchheimer law, with nonlinear function $G$ as in \eqref{G-explicit}. In this case the coefficient $\beta(x)$ in  \eqref{G-explicit} is taken with its relative magnitude  $\Delta\beta/\beta_{min}=1, 10, 100$, where $\Delta\beta$ is the difference between the maximum value of $\beta(x)$ and the minimum value $\beta_{min}$.

The coarse grid is considered to be $20\times20$ where each of the coarse-grid block contains $20\times20$ fine blocks, relative magnitude of the permeability is $\|\Delta k\|/k_{min}=10$, where $\Delta k$ is the difference between the maximum value of $k(x)$ and the minimum value $k_{min}$. For the layered system we consider both fine and coarse grids to be square. For the random system three different cases are considered: $\frac{H_1}{H_2}=\frac{h_1}{h_2}=0.1, 1, 10$. Here $H_1, H_2$ and $h_1, h_2$ are the size of coarse and fine cells correspondingly.

Both analytic averaging formulas are computationally cheap. 
They show different
performances: formulas \eqref{||ks}, \eqref{||Gs}, derived for  flow parallel to layers, are consistently better than formula \eqref{|ks-gs} derived for  flow perpendicular to layers.
From Table~\ref{tab:incomp-|flow} it can be seen that the accuracy of both formulas decreases as the relative magnitude of nonlinear coefficient $\beta$ increases.

\begin{table*}[ht]
\centering
\begin{tabular}{c| c| c | c |c}
\hline\hline
     & $\beta=0$ & $\frac{\Delta\beta}{\beta_{min}}=1$& $\frac{\Delta\beta}{\beta_{min}}=10$ 
& $\frac{\Delta\beta}{\beta_{min}}=100$\\
\hline
Av $\Rightarrow$ & 0 & 0 & 0 & 0 \\ 
\hline
Av $\Uparrow$ & 0 & 5.36e-4 & 2.36e-3 & 1.7e-2 \\ 
\hline
Num  & 0 & 0 & 0 & 0  \\ 
\hline
\end{tabular}
\smallskip
\caption{Upscaling errors for numerical and analytical methods, permeability of Fig.~\ref{fig:khor}}
\label{tab:incomp-||flow}
\end{table*}
\begin{table*}[ht]
\centering
\begin{tabular}{c| c| c | c |c}
\hline\hline
     & $\beta=0$ & $\frac{\Delta\beta}{\beta_{min}}=1$& $\frac{\Delta\beta}{\beta_{min}}=10$ 
& $\frac{\Delta\beta}{\beta_{min}}=100$\\
\hline
Av $\Rightarrow$ & 0 & 8.38e-4 & 9.56e-4 & 7.6e-4 \\ 
\hline
Av $\Uparrow$ & 0 & 0 & 0 & 0 \\ 
\hline
Num  & 0 & 0 & 0 & 0  \\ 
\hline
\end{tabular}
\smallskip
\caption{Upscaling errors for numerical and analytical methods, permeability of Fig.~\ref{fig:kver}}
\label{tab:incomp-|flow}
\end{table*}
\begin{table*}[ht]
\centering
\begin{tabular}{c| c| c | c |c}
\hline\hline
$\frac{H_1}{H_2}=\frac{h_1}{h_2}=10$ & $\beta=0$ & $\frac{\Delta\beta}{\beta_{min}}=1$& $\frac{\Delta\beta}{\beta_{min}}=10$ 
& $\frac{\Delta\beta}{\beta_{min}}=100$\\
\hline
Av $\Rightarrow$ & 8.5e-3 & 7.39e-3 & 8.63e-3 & 2.43e-2 \\ 
\hline
Av $\Uparrow$ & 8.5e-3 & 1.87e-2 & 8.44e-2 & 0.18e-1 \\ 
\hline
Num  & 8.5e-3 & 6.52e-3 & 9.65e-3 & 1.70e-2 \\ 
\hline\hline
$\frac{H_1}{H_2}=\frac{h_1}{h_2}=1$ & $\beta=0$ & $\frac{\Delta\beta}{\beta_{min}}=1$& $\frac{\Delta\beta}{\beta_{min}}=10$ 
& $\frac{\Delta\beta}{\beta_{min}}=100$\\
\hline
Av $\Rightarrow$ & 1.38e-2 & 1.52e-2 & 2.76e-2 & 9.81e-2 \\ 
\hline
Av $\Uparrow$ & 1.38e-2 & 1.35e-2 & 6.55e-2 & 1.20e-1 \\ 
\hline
Num  & 1.38e-2 & 9.53e-3 & 1.25e-2 & 1.60e-2 \\ 
\hline\hline
$\frac{H_1}{H_2}=\frac{h_1}{h_2}=0.1$ & $\beta=0$ & $\frac{\Delta\beta}{\beta_{min}}=1$& $\frac{\Delta\beta}{\beta_{min}}=10$ 
& $\frac{\Delta\beta}{\beta_{min}}=100$\\
\hline
Av $\Rightarrow$ & 1.86e-2 & 2.21e-2 & 3.87e-2 & 1.32e-1 \\ 
\hline
Av $\Uparrow$ & 1.86e-2 & 1.11e-2 & 5.59e-2 & 9.33e-2 \\ 
\hline
Num  & 1.86e-2 & 1.09e-3 & 1.54e-2 & 1.99e-2 \\
\hline
\end{tabular}
\smallskip
\caption{Upscaling errors for numerical and analytical methods, randomly distributed permeability}
\label{tab:incomp-rand}
\end{table*}

\subsection{Numerical Results for Slightly Compressible Fluid}\label{sec:num-compres}
Numerical results for upscaling in case of slightly-compressible flow are presented in Tables \ref{tab:compr-hor}-\ref{tab:compr-randcorr}.
The coarse grid is taken to be $4\times4$ where each of them consists of $64\times64$ fine blocks. 

On the fine scale the pressure is subject to equation \eqref{eq-compr-fine} in the region $\Omega=[0,L_1]\times[0,L_2]$ with the boundary 
conditions
The following boundary conditions are imposed 
\begin{equation*}
 p(L_1,x_2)=0;\quad \tfrac{\partial p}{\partial x_1}(0,x_2)=\tfrac{\partial p}{\partial x_2}(x_1,0)=\tfrac{\partial p}{\partial x_2}(x_1,L_2)=0.
\end{equation*}
We report  the relative errors in the average velocity and the PI between exact and upscaled solution, given by \eqref{vel-err} and 
$|PI-PI^*|/PI$.

The calculations are performed for linear Darcy case and two-term Forchheimer law, with $G(x_1,x_2)$ as in \eqref{G-explicit}.
Four distributions of the fine permeability field $k$ are considered,  see Fig.~\ref{fig:k-distr}.
The commonly used empirical formulas to relate porosity $\phi$ and permeability $k$  are of the form  
\begin{equation*}
 \phi\sim k^{\alpha_0}
\end{equation*}
where $\alpha_0=0.2222$ (Willie and Rose equation),  $\alpha_0=0.2272$ (Timur equation), 
$\alpha_0=1/3$ (Morris and Biggs equation), 
 $\alpha_0=0.1961$ (Berg equation). 
In our calculations  we take $\phi=0.1\cdot k^{\alpha}$ where $\alpha=0.33, 0.25, 0.2$.
Following \cite{AIVW09}  the Forchheimer coefficient $\beta$  is taken
\begin{equation*}
 \beta=\frac{\phi}{k^{1/2}}
\end{equation*}

Numerical results show that the proposed upscaling algorithm provides small 
errors for the upscaled average velocity and productivity index.
In particular, the relative errors are less than 5 \% in all cases.  Though,
we use a different expression for the case of compressible flow compared
to incompressible flow, we observe that 
the velocity errors become larger for nonlinear flows.

\begin{table*}[ht]
\centering
\begin{tabular}{c| c| c ||c |c||c|c}
\hline\hline
&\multicolumn{2}{c||}{$\alpha=1/3$} &\multicolumn{2}{c||}{$\alpha=1/4$}&\multicolumn{2}{c}{ $\alpha=1/5$}\\
\hline
                      & PI err & Vel err&  PI err & Vel err&  PI err&Vel Err \\

\hline
Darcy  & 2.93e-2 & 5.5e-3 & 2.86e-2 & 4.3e-3 & 2.81e-2  & 3.5e-3 \\
\hline\hline
2-Forch & 2.4e-3 & 7.6e-3 & 3.5e-3 & 7.1e-3 & 4.3e-3  & 7.2e-3 \\
\hline
\end{tabular}
\smallskip
\caption{Upscaling errors, permeability of Fig.~\ref{fig:kver}}
\label{tab:compr-ver}
\end{table*}
\vspace{-0.7cm}
\begin{table*}[ht]
\centering
\begin{tabular}{c| c| c ||c |c||c|c}
\hline\hline
&\multicolumn{2}{c||}{$\alpha=1/3$} &\multicolumn{2}{c||}{$\alpha=1/4$}&\multicolumn{2}{c}{ $\alpha=1/5$}\\
\hline
                      & PI err & Vel err&  PI err & Vel err&  PI err&Vel Err \\

\hline
Darcy  & 1.3e-3 & 6.4e-3 & 1.5e-3 & 7.4e-3 & 1.6e-3  & 8.0e-3 \\
\hline\hline
2-Forch & 3.9e-3 & 8.3e-3 & 4.4e-3 & 9.2e-3 & 4.7e-3  & 9.9e-3 \\
\hline
\end{tabular}
\smallskip
\caption{Upscaling errors, permeability of Fig.~\ref{fig:khor}}
\label{tab:compr-hor}
\end{table*}
\vspace{-0.7cm}
\begin{table*}[ht]
\centering
\begin{tabular}{c| c| c ||c |c||c|c}
\hline\hline
&\multicolumn{2}{c||}{$\alpha=1/3$} &\multicolumn{2}{c||}{$\alpha=1/4$}&\multicolumn{2}{c}{ $\alpha=1/5$}\\
\hline
                      & PI err & Vel err&  PI err & Vel err&  PI err&Vel Err \\

\hline
Darcy  & 3.9e-3  & 3.07e-2 &  3.7e-3  & 3.12e-2 & 3.5e-3  & 3.15e-2 \\
\hline\hline
2-Forch & 1.72е-2 & 3.35е-2 & 1.84e-2  & 3.48e-2 & 1.92e-2  & 3.57e-2 \\
\hline
\end{tabular}
\smallskip
\caption{Upscaling errors, permeability of Fig.~\ref{fig:krand}}
\label{tab:compr-randcorr}
\end{table*}
\vspace{-0.7cm}
\begin{table*}[ht]
\centering
\begin{tabular}{c| c| c ||c |c||c|c}
\hline\hline
&\multicolumn{2}{c||}{$\alpha=1/3$} &\multicolumn{2}{c||}{$\alpha=1/4$}&\multicolumn{2}{c}{ $\alpha=1/5$}\\
\hline
                      & PI err & Vel err&  PI err & Vel err&  PI err&Vel Err \\

\hline
Darcy  & 1.98e-2 & 5.7e-3 & 1.9e-2 & 5.5e-3 & 1.84e-2  & 5.5e-3 \\
\hline\hline
2-Forch & 4.4e-3 & 9.7e-3 & 5.3e-3 & 1.05e-2 & 5.9e-3  & 1.10e-2 \\
\hline
\end{tabular}
\smallskip
\caption{Upscaling errors, permeability of Fig.~\ref{fig:kplane}}
\label{tab:compr-plane}
\end{table*}

%

%


\section{Conclusions}
\begin{itemize}
 \item The developed upscaling algorithm for nonlinear steady state problems can can be effectively used for $p$-Laplacian type equations of the form  \eqref{eq-incomp-fine} and \eqref{eq-compr-fine}  and for variety of heterogeneities in the
 domain of computation.

\item  The coarse scale parameters   $\kbs$, $\gs$ and $\Phis$  are determined so that the volumetric average of velocity of the flow in the domain on fine scale and on coarse scale are close enough.

\item The numerical results show that the proposed method  can be 
used to approximate  the Productivity Index (PI) of the well in the bounded domain on the coarse scale.

\item Analytical upscaling formulas in stratified domain are obtained for the nonlinear case. 

\item In our results for the nonlinear problems,  the upscaled parameters  depend on the range of boundary data. 

\item Our results on asymptotic behavior of fully transient velocity and PI and actual numerical computations justify the usage of the coarse scale parameters $\kbs$, $\gs$ and $\Phis$ obtained for the steady state case in the fully transient problem \eqref{degPar}.
 
\end{itemize}

\setcounter{section}{0}
\renewcommand{\thesection}{\Alph{section}}
\section{Appendix}\label{sec:appendix}
\setcounter{equation}{0}
\numberwithin{equation}{section}
In the Appendix we  present the prospective improvement of the upscaling method, which in some cases can give better results with only small increase in computational costs (see Tables~\ref{tab:compr-hor-full}-\ref{tab:compr-randcorr-full}).  

Here we consider the upscaled $\kbs$ and $\Phis$ to be of a form
\begin{align}
&\kbs=
\begin{bmatrix}
K_{11}\cdot \kbar & K_{12}\\
K_{21} & K_{22}\cdot \kbar 
\end{bmatrix},\label{ks}\\
&\qquad\quad \text{where}\quad \kbar=K_0+K_1(x_1-c_1)+K_2(x_2-c_2);\label{kbar}\\
&\Phis=\Phi_A+\Phi_B(x_1-c_1)+\Phi_C(x_2-c_2)+\Phi_D[(x_1-c_1)^2-(x_2-c_2)^2]. \label{phis}
\end{align}
Here $c=(c_1,c_2)$ is the central point of the coarse cell $\Omega_c$ and $K_{11}, K_{12}, K_{21}, K_{22}$, $K_0, K_1, K_2$, $\Phi_A, \Phi_B, \Phi_C$ and $\Phi_D$ are constants to be determined.

First, the permeability tensor $\kbs$ is determined, and then it is used to upscale the porosity $\phi$.  

I. {\it Permeability $\kbs$.}  We first obtain the polynomial $\kbar$ in \eqref{kbar} 
as the least square approximation of permeability $k(x_1,x_2)$:
%
%
\begin{equation}\label{kbar-coeff}
   K_0=\frac1{|\Omega_c|}\int_{\Omega_c}k(x_1,x_2)\,d\Omega_c;\quad
  K_i=\dfrac{\int_{\Omega_c}k(x_1,x_2)(x_i-c_i)\,d\Omega_c}{\int_{\Omega_c}(x_i-c_i)^2\,d\Omega_c},\ \ i=1,2.
\end{equation}

We now will use the constant elements $K_{11}$, $K_{12}$, $K_{21}$ and  $K_{22}$ of matrix $\kbs$ to ``correct'' $\kbar(x_1,x_2)$ so that the average velocities
 on the fine and coarse scales are the same. For this purpose we will modify our approach presented in Sec.~\ref{sec:algorithm}.

We consider two linear fine-scale equations with zero RHS, the exact equation \eqref{eq_p_k} and the averaged equation
\begin{equation}
\nabla\cdot(\kbar(x_1,x_2)\nabla P)=0.\label{eq_P_k}
\end{equation}
Each equation we solve twice.
Namely, let $p_1$ and $P_1$ be the solutions of \eqref{eq_p_k} and \eqref{eq_P_k}, correspondingly, subject to boundary conditions (\ref{k-bc}$_1$) and let
$p_2$ and $P_2$ be the solutions of \eqref{eq_p_k} and \eqref{eq_P_k}, correspondingly, subject to boundary conditions (\ref{k-bc}$_2$).

We now equate the velocity averages on fine and coarse scale, with the coarse scale velocity $-\kbar\nabla P_i$, $i=1,2$, ``corrected'' with the elements of the matrix $\kbs$:
\begin{equation}\label{u=us-for-k}
\frac1{|\Omega_c|}\int_{\Omega_c}-k\nabla p_i\,d\Omega_c= \frac1{|\Omega_c|}\int_{\Omega_c}-\kbs\nabla P_i\,d\Omega_c, \qquad i=1,2.
\end{equation}
Solving the linear system of four equations \eqref{u=us-for-k} gives the values of $K_{11}$, $K_{12}$, $K_{21}$, $K_{22}$.

II. {\it Porosity $\Phis$.} 
To find the coefficients $\Phi_A$, $\Phi_B$, $\Phi_C$ and $\Phi_D$ in \eqref{phis} we
consider two equations:
the upscaled linear equation with function \eqref{phis} in RHS
\begin{equation}\label{eq_av_phis}
 -\nabla\cdot(\kbs\nabla P)=\Phi_A+\Phi_B(x_1-c_1)+\Phi_C(x_2-c_2)+\Phi_D[(x_1-c_1)^2-(x_2-c_2)^2],
\end{equation}
and the exact linear equation 
\begin{equation}\label{bvp_p_P}
 -\nabla\cdot(k(x_1,x_2)\nabla p_P)=\phi(x_1,x_2).
\end{equation}
Zero Dirichlet boundary conditions are imposed:
\begin{equation*}
  P|_{\partial \Omega_c}=p_P|_{\partial \Omega_c}=0\notag.
\end{equation*}

Coefficients $\Phi_A$, $\Phi_B$, $\Phi_C$ and $\Phi_D$ are determined so that the boundary fluxes through the faces of the coarse cell corresponding to the upscaled solution $P$ and 
 fine scale solution $p_P$ are equal.

Due to the linearity of equation  \eqref{eq_av_phis} solution $P$ is the linear combination
\begin{equation}\label{phi-eq-4}
 P=\Phi_A\cdot P_A+\Phi_B\cdot P_B+\Phi_C\cdot P_C+\Phi_D\cdot P_D,
\end{equation}
where $P_{\lambda}$, $\lambda=A,B,C,D$, are the solutions of BVPs
\begin{align}\label{eqs-ABCD}
& -\nabla\cdot(\kbs\nabla P_{\lambda})=f_\lambda(x_1,x_2),\\
&P_{\lambda}|_{\partial \Omega_c}=0\notag
\end{align}
with the RHS
\begin{equation*}
 f_{\lambda}(x_1,x_2)=
\begin{cases}
 1 &\quad\text{for}\quad \lambda=A,\\
  x_1-c_1 &\quad\text{for}\quad \lambda=B,\\
  x_2-c_2 &\quad\text{for}\quad \lambda=C,\\
 (x_1-c_1)^2-(x_2-c_2)^2 &\quad\text{for}\quad \lambda=D.
\end{cases}
\end{equation*}
The boundary fluxes through four faces of the coarse cell are then related by the same expression as \eqref{phi-eq-4}. 
Solving the resulting system of four equations for $\Phi_\lambda$, $\lambda=A,B,C,D$, we obtain the expression for $\Phis$.

The corresponding numerical results are presented in Tables \ref{tab:compr-hor-full}-\ref{tab:compr-randcorr-full}.
For the reader's convenience the results presented in   Tables \ref{tab:compr-hor}-\ref{tab:compr-randcorr} are included here once again for the comparison.
We compare errors for the upscaling algorithm with $\ks$ and $\Phis$ calculated using different approaches described above.
Four cases are considered for each of linear and nonlinear case.
\begin{enumerate}
 \item[(i)]  $\ks$ and $\Phis$ are calculated as in Sec.~\ref{sec:algorithm} (denoted as $\kbs$ - C; $\Phis$ - C);
\item[(ii)] $\ks$ is as in \eqref{ks}, while $\Phis$ is a constant \eqref{phia} (denoted as $\kbs$ - P; $\Phis$ - C);
\item[(iii)] $\ks$ is as in Sec.~\ref{sec:algorithm}, while $\Phis$ is as in \eqref{phis} (denoted as $\kbs$ - C; $\Phis$ - P);
\item[(iv)] $\ks$ is as in \eqref{ks} and $\Phis$ is as in \eqref{phis} (denoted as $\kbs$ - P; $\Phis$ - P).
\end{enumerate}

Approach (iv) provides consistently better results for both error in PI and velocity. 
Though it is computationally more expensive, the increase is negligible in nonlinear case where the main computational expense comes from obtaining the function $\gs$.
Approach (ii)  is routinely comparable to (iv), and for some permeability distributions ($k$ as in Fig.~\ref{fig:kver}, $\alpha=1/4, 1/5$) is even better. With that it only amounts to computation of coefficients \eqref{kbar-coeff}.
As expected,  approach (iii) is routinely the worst, as there is not sufficient information for correct estimation of coefficients of polynomial \eqref{phis}.

\section*{Acknowledgments}
The authors are thankful to Dr. Luan Hoang for his valuable discussions, 
suggestions and recommendations. 
The research of this paper was supported by the NSF grant DMS-0908177.

\bibliographystyle{elsarticle-num}
\bibliography{reference}

\begin{table*}[ht]
\centering
\begin{tabular}{c| c| c ||c |c||c|c}
\hline\hline
&\multicolumn{2}{c||}{$\alpha=1/3$} &\multicolumn{2}{c||}{$\alpha=1/4$}&\multicolumn{2}{c}{ $\alpha=1/5$}\\
\hline
                      & PI err & Vel err&  PI err & Vel err&  PI err&Vel Err \\

\hline
Lin. $\kbs$ - P; $\Phis$ - P & 5.5e-3 & 1.9e-3 & 5.6e-3 & 2.1e-3 & 5.6e-3  & 2.2e-3 \\
Lin. $\kbs$ - P; $\Phis$ - C & 9.3e-3 & 5.5e-3 & 8.6e-3 & 4.3e-3 & 8.2e-3  & 3.5e-3 \\
Lin. $\kbs$ - C; $\Phis$ - P & 2.17e-2 & 7.9e-3 & 2.17e-2 & 8.1e-3 & 2.18e-2  & 8.1e-3 \\
Lin. $\kbs$ - C; $\Phis$ - C & 2.93e-2 & 5.5e-3 & 2.86e-2 & 4.3e-3 & 2.81e-2  & 3.5e-3 \\
\hline\hline
$\kbs$ - P; $\Phis$ - P & 1.6e-3 & 6.7e-3 & 2.5e-3 & 7.2e-3 & 3.3e-3  & 7.6e-3 \\
$\kbs$ - P; $\Phis$ - C & 1.4e-3 & 7.2e-3 & 2.3e-3 & 6.6e-3 & 3.0e-3  & 6.5e-3 \\
$\kbs$ - C; $\Phis$ - P & 2.4e-3 & 1.08e-2 & 3.5e-3 & 1.13e-2 & 4.4e-3  & 1.17e-2 \\
$\kbs$ - C; $\Phis$ - C & 2.4e-3 & 7.6e-3 & 3.5e-3 & 7.1e-3 & 4.3e-3  & 7.2e-3 \\
\hline
\end{tabular}
\smallskip
\caption{Comparison of upscaling errors, permeability of Fig.~\ref{fig:kver}}
\label{tab:compr-ver-full}
\end{table*}
\vspace{-0.7cm}
\begin{table*}[ht]
\centering
\begin{tabular}{c| c| c ||c |c||c|c}
\hline\hline
&\multicolumn{2}{c||}{$\alpha=1/3$} &\multicolumn{2}{c||}{$\alpha=1/4$}&\multicolumn{2}{c}{ $\alpha=1/5$}\\
\hline
                      & PI err & Vel err&  PI err & Vel err&  PI err&Vel Err \\

\hline
Lin. $\kbs$ - P; $\Phis$ - P & 1.8e-3 & 2.1e-3 & 1.8e-3 & 2.4e-3 & 1.7e-3  & 2.5e-3 \\
Lin. $\kbs$ - P; $\Phis$ - C & 2.5e-4 & 3.9e-3 & 1.4e-4 & 2.9e-3 & 6.6e-5  & 2.3e-3 \\
Lin. $\kbs$ - C; $\Phis$ - P & 3.6e-3 & 1.68e-2 & 3.6e-3 & 1.70e-2 & 3.7e-3  & 1.72e-2 \\
Lin. $\kbs$ - C; $\Phis$ - C & 1.3e-3 & 6.4e-3 & 1.5e-3 & 7.4e-3 & 1.6e-3  & 8.0e-3 \\
\hline\hline
$\kbs$ - P; $\Phis$ - P & 3.5e-3 & 6.8e-3 & 4.0e-3 & 7.2e-3 & 4.4e-3  & 7.4e-3 \\
$\kbs$ - P; $\Phis$ - C & 3.7e-3 & 7.3e-3 & 4.2e-3 & 7.1e-3 & 4.5e-3  & 7.0e-3 \\
$\kbs$ - C; $\Phis$ - P & 3.8e-3 & 1.76e-2 & 4.3e-3 & 1.80e-2 & 4.6e-3  & 1.82e-2 \\
$\kbs$ - C; $\Phis$ - C & 3.9e-3 & 8.3e-3 & 4.4e-3 & 9.2e-3 & 4.7e-3  & 9.9e-3 \\
\hline
\end{tabular}
\smallskip
\caption{Comparison of upscaling errors, permeability of Fig.~\ref{fig:khor}}
\label{tab:compr-hor-full}
\end{table*}
\vspace{-0.7cm}
\begin{table*}[ht]
\centering
\begin{tabular}{c| c| c ||c |c||c|c}
\hline\hline
&\multicolumn{2}{c||}{$\alpha=1/3$} &\multicolumn{2}{c||}{$\alpha=1/4$}&\multicolumn{2}{c}{ $\alpha=1/5$}\\
\hline
                      & PI err & Vel err&  PI err & Vel err&  PI err&Vel Err \\

\hline
Lin. $\kbs$ - C; $\Phis$ - C & 3.9e-3  & 3.07e-2 &  3.7e-3  & 3.12e-2 & 3.5e-3  & 3.15e-2 \\
Lin. $\kbs$ - P; $\Phis$ - C & 5.9е-3 & 1.27e-2 & 6.1е-3 & 1.23e-2 & 6.1е-3  & 1.20e-2 \\
Lin. $\kbs$ - C; $\Phis$ - P & 8.4е-4 & 3.88e-2 & 7.7е-4 & 3.89e-2 & 7.2е-4  & 3.90e-2 \\
Lin. $\kbs$ - P; $\Phis$ - P & 5.9е-3 & 1.23e-2 & 5.9е-3 & 1.23e-2 & 6.0е-3  & 1.23e-2 \\
\hline\hline
$\kbs$ - C; $\Phis$ - C & 1.72е-2 & 3.35е-2 & 1.84e-2  & 3.48e-2 & 1.92e-2  & 3.57e-2 \\
$\kbs$ - P; $\Phis$ - C & 1.10e-2 & 1.36e-2 & 1.26e-2 & 1.36e-2 & 1.38e-2  & 1.4e-2 \\
$\kbs$ - C; $\Phis$ - P & 9.6e-3 & 3.9e-2  & 1.11e-2 & 3.9e-2 & 1.21e-2 & 3.9e-2  \\
$\kbs$ - P; $\Phis$ - P & 1.49e-4 & 2.0e-2 & 1.63e-2 & 2.12e-2 &  1.39e-2 & 1.43e-2 \\
\hline
\end{tabular}
\smallskip
\caption{Comparison of upscaling errors, permeability of Fig.~\ref{fig:krand}}
\label{tab:compr-randcorr-full}
\end{table*}
\vspace{-0.7cm}
\begin{table*}[ht]
\centering
\begin{tabular}{c| c| c ||c |c||c|c}
\hline\hline
&\multicolumn{2}{c||}{$\alpha=1/3$} &\multicolumn{2}{c||}{$\alpha=1/4$}&\multicolumn{2}{c}{ $\alpha=1/5$}\\
\hline
                      & PI err & Vel err&  PI err & Vel err&  PI err&Vel Err \\

\hline
Lin. $\kbs$ - P; $\Phis$ - P & 6.2e-5 & 4.5e-5 & 6.7e-5 & 4.7e-5 & 7.2e-5  & 4.8e-5 \\
Lin. $\kbs$ - P; $\Phis$ - C & 3.1e-3 & 4.5e-3 & 2.3e-3 & 3.4e-3 & 1.8e-3  & 2.7e-3 \\
Lin. $\kbs$ - C; $\Phis$ - P & 1.21e-3 & 9.4e-3 & 1.21e-2 & 9.5e-3 & 1.2e-2  & 9.6e-3 \\
Lin. $\kbs$ - C; $\Phis$ - C & 1.98e-2 & 5.7e-3 & 1.9e-2 & 5.5e-3 & 1.84e-2  & 5.5e-3 \\
\hline\hline
$\kbs$ - P; $\Phis$ - P & 3.5e-3 & 5.9e-3 & 4.3e-3 & 6.4e-3 & 4.9e-3  & 6.8e-3 \\
$\kbs$ - P; $\Phis$ - C & 3.7e-3 & 6.9e-3 & 4.5e-3 & 6.6e-3 & 5.0e-3  & 6.7e-3 \\
$\kbs$ - C; $\Phis$ - P & 4.0e-3 & 1.32e-2 & 4.9e-3 & 1.38e-2 & 5.6e-3  & 1.42e-2 \\
$\kbs$ - C; $\Phis$ - C & 4.4e-3 & 9.7e-3 & 5.3e-3 & 1.05e-2 & 5.9e-3  & 1.10e-2 \\
\hline
\end{tabular}
\smallskip
\caption{Comparison of upscaling errors, permeability of Fig.~\ref{fig:kplane}}
\label{tab:compr-plane-full}
\end{table*}

\end{document}